\documentclass{amsart}
\usepackage{amsmath, amsthm, amssymb, amscd, epsfig, mathrsfs, eucal, palatino}
\input rlepsf.sty

\DeclareMathOperator\height{height}

\begin{document}

\newtheorem{theorem}{Theorem}[section]
\newtheorem{lemma}[theorem]{Lemma}
\newtheorem{corollary}[theorem]{Corollary}
\newtheorem{conjecture}[theorem]{Conjecture}
\newtheorem{question}[theorem]{Question}
\newtheorem{problem}[theorem]{Problem}
\newtheorem*{claim}{Claim}
\newtheorem*{criterion}{Criterion}
\newtheorem*{main_thm}{Theorem A}
\newtheorem*{precise_main_thm}{Theorem A'}

\theoremstyle{definition}
\newtheorem{definition}[theorem]{Definition}
\newtheorem{construction}[theorem]{Construction}
\newtheorem{notation}[theorem]{Notation}

\theoremstyle{remark}
\newtheorem{remark}[theorem]{Remark}
\newtheorem{example}[theorem]{Example}

\def\id{\text{id}}
\def\Z{\mathbb Z}
\def\R{\mathbb R}
\def\Q{\mathbb Q}
\def\fix{\textnormal{fix}}
\def\PL{\textnormal{PL}}
\def\homeo{\textnormal{Homeo}}
\def\diffeo{\textnormal{Diffeo}}
\def\rot{\textnormal{rot}}
\def\inte{\textnormal{int}}

\title{Denominator bounds in Thompson-like groups and flows}
\author{Danny Calegari}

\begin{abstract}
Let $T$ denote Thompson's group of piecewise $2$-adic linear homeomorphisms of
the circle. Ghys and Sergiescu showed that the
rotation number of every element of $T$ is rational, but their proof is
very indirect. We give here a short, direct proof using train tracks,
which generalizes to
elements of $\PL^+(S^1)$ with rational break points and derivatives 
which are powers of some fixed integer, and also to certain
flows on surfaces which we call {\em Thompson-like}. We also obtain an explicit
upper bound on the smallest period of a fixed point in terms of data
which can be read off from the combinatorics of the homeomorphism.
\end{abstract}

\address{Department of Mathematics \\ California Institute of Technology \\
Pasadena CA, 91125}
\email{dannyc@its.caltech.edu}
\date{1/25/2007, Version 0.08}

\maketitle

\section{Introduction}

In \cite{Ghys_Sergiescu}, Ghys and Sergiescu studied Thompson's group
$T$ of homeomorphisms of the circle from a number of points of view. This
group was introduced in unpublished notes by Thompson, and is defined 
as the subgroup of $\homeo^+(S^1)$ consisting of homeomorphisms taking dyadic rationals
to dyadic rationals which are piecewise linear, where the break points are dyadic rationals
(i.e. numbers of the form $p2^q$ for $p,q \in \Z$), and where the derivatives
are all of the form $2^q$ for $q \in \Z$. 

One of the main theorems in \cite{Ghys_Sergiescu} is that the {\em rotation
number} of every element of $T$ is rational. Recall Poincar\'e's definition \cite{Poincare}
of {\em rotation number} for an element $g \in \homeo^+(S^1)$. Let
$\tilde{g}$ be any lift of $g$ to $\homeo^+(\R)$, and define
$$\rot(\tilde{g}) = \lim_{n \to \infty} \frac {\tilde{g}^n(0)} n$$
Then $\rot(g) = \rot(\tilde{g}) \pmod \Z$; i.e. $\rot(g) \in S^1$.
Different lifts $\tilde{g}_1,\tilde{g}_2$ of $g$ satisfy 
$$\rot(\tilde{g}_1) - \rot(\tilde{g}_2) \in \Z$$
so this is well-defined.

In \cite{Ghys_circles}, Ghys says ``the proof (of rationality)
is very indirect and there is a need for a better proof". One such argument was given
by I. Liousse \cite{Liousse1},\cite{Liousse2}. We also recently learned that
V. Kleptsyn has an approach to understanding rationality in Thompson's group
using automata, which is distinct from, but not unrelated to, the approach
in this paper; see \cite{Kleptsyn}.

The argument of Ghys--Sergiescu
is a proof by contradiction: they show there is a morphism $\rho:T \to \diffeo^\infty(S^1)$
which is semi-conjugate to the natural (topological) action, and which has an
exceptional minimal set. Rotation number is invariant under semi-conjugacy.
Therefore the existence of an element of $T$ with irrational
rotation number would contradict Denjoy's theorem (that every $C^2$ diffeomorphism
of $S^1$ with an irrational rotation number has dense orbits). Liousse's proof is more
straightforward, but is still nonconstructive, and is still a proof by contradiction
(to Denjoy's inequalities).

\vskip 12pt

It is well-known and easy to show (see e.g. \cite{Poincare}) that an element
$g \in \homeo^+(S^1)$ has a periodic point of (least) period $q$ if and only if
its rotation number is $p/q$ for some coprime pair of integers $p,q$. In this note,
we supply a direct proof of rationality of rotation number for certain
$\PL$ homeomorphisms by directly finding a periodic point. In fact, our
argument applies more generally than the argument of \cite{Ghys_Sergiescu}, though
not more generally than the argument of \cite{Liousse2}. On the other hand, since it
is constructive, we obtain explicit bounds on the denominator of
the rotation number in terms of the combinatorics of the original map, which are
not obtained in either \cite{Ghys_Sergiescu} or \cite{Liousse2}.

\vskip 12pt

The main rationality theorem, proved in \S~2, is as follows:

\begin{main_thm}
Let $t$ be an element of $\PL^+(S^1)$ mapping rationals to rationals,
with rational break points, and
with derivatives of the form $n^q$ for $q \in \Z$ and some fixed $n$.
Then the rotation number of $t$ is rational.
\end{main_thm}

This theorem is also proved in \cite{Liousse2}.
We remark that {\em some} hypothesis on $t$ is necessary, since there are examples
of elements of $\PL^+(S^1)$ mapping rationals to rationals,
with rational break points and rational derivatives,
whose rotation numbers are irrational. 
However, it is possible that the conditions on $t$ could still be relaxed further
(c.f. Question~\ref{rational_derivative_question}). 
See \cite{Liousse2} and also \cite{Herman} or \cite{Boshernitzan}.

\vskip 12pt

Given $t \in \PL^+(S^1)$ as in the statement of the theorem, 
we define the {\em height} of $t$ as follows. Let $m$ be the least integer
such that a subdivision of $S^1$ into $m$ intervals of the form $\left[ \frac p m ,\frac {p+1} m\right]$ 
is {\em Markov} for $t$; i.e. $t$ either takes a strip of $n^k$ consecutive 
intervals linearly to a single interval
by contraction, or takes a single interval and stretches it linearly over $n^k$ consecutive
intervals by expansion, for various $k$. Then define $\height(t) = m$.

\begin{remark}
Note that $m$ as we have defined it is bounded by the least common multiple of the
denominators of the break points of $t$ and their images (assuming there is at least one break
point). This remark justifies the existence of $m$
and gives a direct estimate for its size.
\end{remark}

In \S~3 we obtain a straightforward bound on denominator of rotation number as follows:

\begin{precise_main_thm}
Let $t$ be an element of $\PL^+(S^1)$ mapping rationals to rationals, 
with rational break points, and with
derivatives of the form $n^q$ for $q \in \Z$ and some fixed $n$.
Suppose we have
$$\height(t) = m$$
Then $t$ has a periodic point of period at most $n^m\cdot m$.
\end{precise_main_thm}

Of course, Theorem A' implies Theorem A.

\vskip 12pt

Finally, in \S~4 we construct some examples with periodic points with long periods,
complementing the estimate in Theorem A'.

\subsection{Acknowledgements}
While writing this paper, I was partially supported by a Sloan Research Fellowship, and
NSF grant DMS-0405491. I'm grateful for comments from \'Etienne Ghys and Collin Bleak, and
for some substantial corrections by the anonymous referee.
I am also grateful to Isabelle Liousse for her comments on
an early version of this paper, and to her and Victor Kleptsyn 
for forwarding me their relevant preprints.

\section{Train tracks}

To prove Theorem A, it suffices to show that $t$ as in the statement of the
theorem has a periodic orbit.

Throughout this section, for concreteness and ease of exposition,
we will concentrate on the case $n=2$. This includes (but is more general
than) the case of the classical Thompson group, although
our argument goes through essentially verbatim for general $n$.

\vskip 12pt

We fix $t$ as in the statement of the theorem with rational break points and
all derivatives powers of $2$. Here is a summary of the proof.
We associate to $t$ a {\em train track} which carries its dynamics; by analyzing the
combinatorics of the train track we show that we can either {\em split off a circle}
in which case some power of $t$ is equal to the identity on some segment, or
else we can find an {\em attracting cycle} in which case some (possibly negative)
power of $t$ has a periodic orbit which is attracting on at least one side. This
will complete the proof.

\vskip 12pt

Train tracks are introduced in \cite{Thurston_notes} as a combinatorial tool for studying
one dimensional dynamics on surfaces (similar objects were introduced earlier by Dehn and
Nielsen). A train track $\tau$ is a graph with a $C^1$ combing
at every vertex which comes with an embedding in a surface. The vertices of
$\tau$ are called the {\em switches}. We now show how to associate a train track $\tau$
in a torus to our element $t$.

\vskip 12pt

By the defining properties of $t$, there is a least integer $m$ such that there is
a Markov partition for $t$ of the simple form
$$S^1 = I_1 \cup I_2 \cup \dots \cup I_m$$
where each $I_i$ has length $1/m$. The element $t$ acts in two
ways: by taking a strip of $2^k$ consecutive intervals and mapping them linearly
to a single interval (contraction), 
or by taking a single interval and stretching it out linearly
over $2^k$ consecutive intervals (expansion), for various positive integers $k$.
Recall that we are calling $m$ the {\em height} of $t$.

The mapping torus of $t$ is literally a (two dimensional)
torus which we denote by $F$. In $F$ we
construct an oriented train track $\tau$ by gluing {\em intervals}. We take one 
oriented interval $e_i$ for each $I_i$, and we glue the intervals together at their
endpoints in a pattern determined by the dynamics of $t$. As a convention, we
think of the circle $S^1$ as being embedded ``horizontally" in $F$, and the
edges $e_i$ as being embedded ``vertically". If $t(I_i) \cap I_j$ has
nonempty interior, then we identify the positive endpoint of $e_i$ to the negative 
endpoint of $e_j$. The orientation on the $e_i$ determines the combing at the
vertices. In this way we obtain an orientable train track $\tau$, 
which comes with a natural embedding in $F$. We emphasize that
a single edge of $\tau$ might consist of many intervals. We do not need to
keep track of intervals in this section, but they will be important in \S~3 when
we try to estimate periods of periodic points.

\begin{figure}[ht]
\center{\scalebox{.3}{\includegraphics{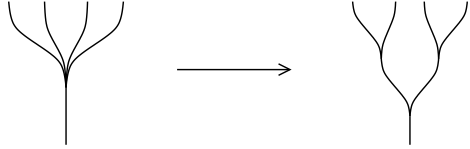}}}
\caption{Split open a $5$-valent switch to three $3$-valent switches}\label{split1}
\end{figure}

At each switch of $\tau$
we either have $2^k$ incoming edges and one outgoing edge for each contraction,
or one incoming edge and $2^k$ outgoing edges for each expansion. Notice that
$k$ may vary from switch to switch. 
We remedy this in the following way: a switch
with $2^k+1$ incident edges can be split open locally to $2^k-1$ switches, each with $3$
incident edges. We say that we are {\em resolving} the (high valence) switches by
this process; see Fig.~\ref{split1} for an example.

We refer to the resolved train track as $\tau_0$. Note that every switch of $\tau_0$
is $3$-valent. We remark that $\tau_0$ could well be disconnected (for that matter,
$\tau$ itself could be disconnected).
Notice that $\tau_0$ is contained in $F$ in such a way that the edges are all transverse
to a foliation of $F$ by meridians, so that an oriented edge of $\tau_0$ points in
a well-defined direction around $F$.
Notice too that the (unparameterized)
dynamics of $t$ can be recovered completely from the combinatorics
of $\tau_0$: we associate to each edge a Euclidean rectangle of width $1$ foliated
by vertical lines. At each switch we attach the (foliated) mapping cylinder
of the linear map $[0,1] \to [0,2]$ to the horizontal boundaries of the edges. 
This gives a foliated surface with boundary
which comes with an embedding in $F$, and which can also be arranged so that
leaves are transverse to the foliation of $F$ by meridians; 
by collapsing complementary regions we get
precisely the (foliated) mapping torus of the homeomorphism $t$; call this operation the
{\em realization} of the combinatorial train track $\tau_0$.

\vskip 12pt

\begin{figure}[ht]
\center{\scalebox{.3}{\includegraphics{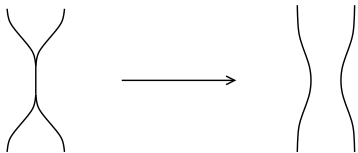}}}
\caption{Sinks can be split open to give a new train track with the same realization}\label{split2}
\end{figure}

Now, at each switch $v$ there is a well defined {\em contracting direction} which
points along the edge which is isolated on its side of $v$. Let $e(v)$ be this
(directed) edge, and let $a(v)$ be the other endpoint of $e(v)$. So $e(v)$ points from
$v$ to $a(v)$. Suppose there is some $v$ such that $e(a(v))$ is the same underlying
edge as $e(v)$, but with the opposite orientation. In this case, if $e$ denotes
the underlying (undirected) edge, we call $e$ a {\em sink}. The key point is
that sinks can be split open to give a new train track whose realization still recovers
the dynamics of $t$. There is only one (obvious) way to do this; see Fig.~\ref{split2}
for an illustration.

Let $\tau_i$ be the result of splitting $\tau_{i-1}$ open along a sink, if one
exists. We split open all sinks until there are none left. Each splitting reduces the number
of vertices by two, so this process must eventually terminate. We denote the
sinkless train track we ultimately obtain by $\tau'$. Note that the realization
of every $\tau_i$, and therefore of $\tau'$, still recovers the dynamics of $t$.
After the splitting, some component of $\tau'$ might be a circle with no switches; 
in such a case we say we have {\em split off a circle}. 
Evidently the realization of a circle is foliated
by periodic orbits. This happens exactly when some power of $t$ fixes a nonempty
interval in $S^1$.

Otherwise, there is still some switch $v$. Since $\tau'$ has only finitely many
switches, the sequence $v,a(v),a^2(v) \dots$ is eventually periodic. Since $\tau'$
is sinkless, the oriented edges $e(a^i(v))$ all point in the same direction around
$F$. It follows that if $v,a(v),\dots,a^n(v)=v$ is a periodic sequence, the union
$$\gamma(v) := e(v) \cup e(a(v)) \cup \dots \cup e(a^{n-1}(v))$$
is an embedded circle in $\tau'$, oriented coherently by the orientation on each edge,
and always pointing in the same direction around $F$.

We call $\gamma(v)$ an {\em attracting cycle}. Notice that all the way around the
realization of an attracting cycle, the linear map at the switches is {\em contracting},
and therefore the realization of an attracting cycle contains a periodic orbit which
is contracting on at least one side. Depending on the orientation of $\gamma(v)$,
this corresponds to a periodic orbit which is attracting on at least one side either
for $t$ or for $t^{-1}$. It follows that $t$ has a periodic orbit, and therefore
the rotation number $\rot(t)$ is rational.

\vskip 12pt

If we replace $2$ by $n$ above, then switches in $\tau$ are $n^k+1$ valent, and
each such switch can be resolved to a union of $\frac {n^k - 1} {n-1}$ switches,
each of valence $n+1$, to produce $\tau_0$. The realization of $\tau_0$ is
obtained by gluing foliated Euclidean rectangles associated to each edge
by attaching the mapping cylinder of the linear map $[0,1] \to [0,n]$ at each
switch. The only point that needs stressing is that sinks of $\tau_{i-1}$
can still be split open to produce $\tau_i$, 
since the two attaching maps in the realization associated to the endpoints of a sink are
still inverse to each other. This proves Theorem A in general.

\begin{remark}
If a general element $t \in \PL^+(S^1)$ has
derivatives and break points which are rational, we can still associate a train track to
$t$ where at every switch there are $p$ incoming and $1$ outgoing edge, and
the realization is obtained by gluing the mapping cylinder of the linear map
$[0,1] \to [0,p]$, for {\em varying} $p$. Note that multiplication by $p/q$ can be realized
as the composition of multiplication by $p$, and division by $q$.
The problem is that one might have a sink where the valences are
{\em different} at the two endpoints. There is no obvious combinatorical simplification
which can be done in such a case, even when the derivatives are all of the
form $r^k$ for some fixed $r$, where $r = p/q$ is rational but not integral 
(see Question~\ref{rational_derivative_question}).
\end{remark}

\begin{remark}
Let $S$ be a closed orientable surface, and let $\tau \subset S$ be an orientable
train track with every switch of valence $3$. The realization of $\tau$ gives
a possibly singular foliation on $S$ with orientable leaves, which defines
an (unparameterized) flow on $S$ that we say is {\em Thompson-like}.
More generally, we call such a flow Thompson-like 
if every switch has valence $n+1$ with $n$ edges on one side, for some fixed $n$.
The argument in this section shows that every orbit in a Thompson-like flow is
either periodic, or accumulates on a periodic orbit.
\end{remark}

\section{Bounding denominators}

To prove the stronger Theorem A' we must analyze the complexity of a split off
circle or an attracting cycle. To do this, we must relate the complexity
of the sinkless train track obtained by the argument of \S~2 to the original
train track.

\vskip 12pt

Let $\tau$ denote the original train track. We distinguish between {\em edges}
of $\tau$ and {\em intervals} which correspond to the original $e_i$, and which
each wrap exactly once around the mapping torus. By definition, an edge of
the train track has both vertices at switches; a single edge may be composed of
many intervals. The period of the periodic
cycle we finally identify will be equal to the number of intervals that it contains.

The train track $\tau$ has exactly $m$ intervals, and has switches of
valence $n^k + 1$ for various $k$. 

We resolve a switch of valence $n^k+1$ to
$\frac {n^k - 1} {n-1}$ switches of valence $n+1$, thereby creating
$\frac {n^k - n} {n-1}$ new edges. We refer to these new edges as {\em infinitesimal
edges}; since the resolution is performed locally, we assume the infinitesimal edges
are as short as we like, and no consecutive sequence of them is long enough to wrap
around the mapping torus. So we can still measure the period of a periodic cycle
by counting the number of intervals it contains, and ignoring the infinitesimal
edges. Observe that after this resolution we obtain a train track which we
call $\tau_0$ with
every switch of valence $n$, in which there are at most $2m$ switches, and
exactly $m$ intervals.

Splitting open a sink produces a new train track with $n-1$ new edges and $2$ fewer
switches. Each edge is a concatenation of intervals and infinitesimal edges, each
of which is replaced by $n-1$ parallel copies after splitting open. If
we let $\tau_i$ be obtained from $\tau_{i-1}$ by splitting open a single sink,
then if there are $m_{i-1}$ intervals in $\tau_{i-1}$, the edge which is
split open contains at most $m_{i-1}$ intervals, and there will be at most
$n\cdot m_{i-1}$ intervals in $\tau_i$. Since $m_0 = m$, and we split open a 
total of at most $m$ sinks, it follows that
if $\tau'$ is the ultimate sinkless train track, the number of intervals in
$\tau'$ is at most $n^m\cdot m$.

Since $\tau'$ is sinkless, it contains an embedded circle or attracting cycle,
and we are done.

\section{Examples}\label{example_section}

In this section we give some examples. Note that the arguments in \S~3 do not
use the embedding of $\tau$ in a torus $F$. Some of the examples in this section
can be realized in a torus, and some cannot.

\begin{example}
The map $t$ might be a rotation of order $m$. The train track $\tau$ is a circle
made up of $m$ intervals.
\end{example}

\begin{example}
Fix some $k$. Let $m = 2^k+s+2$ and define $t$ to be linear on the following intervals:
$$t:\left[\frac 0 m ,\frac {2^k} m \right] \to \left[ \frac {2^k} m, \frac {2^k + 1} m \right]$$
$$t:\left[\frac {2^k} m, \frac {2^k + s} m \right] \to \left[ \frac {2^k+1} m, \frac {2^k + s + 1} m \right]$$
$$t:\left[\frac {2^k+s} m, \frac {2^k+s+1} m\right] \to \left[ \frac {2^k+s+1} m,\frac {2^k+s+1+2^k} m \right]
= \left[ \frac {-1} m, \frac {2^k -1} m \right]$$
$$t:\left[\frac {-1} m, \frac 0 m\right] \to \left[\frac{2^{k}-1}{m}, \frac{2^{k}}{m}\right]$$
The associated train track has two switches of valence $2^k+1$ bounding a sink of length $s+1$
(except when $s=0,k=1$).
The two bushy sides of the switches are glued together with a ``twist". When $\tau$ is completely
split open, the result $\tau'$ is a single circle containing $2^k + 2^k\cdot (s+1) + 1$ intervals, so the
period is $2^k \cdot(s+2)+1$. If we choose $s = 2^k-2$, then $t$ is actually contained in Thompson's
group $T$, and the order of the periodic point is $O(m^2)$.
\end{example}

\begin{example}\label{very_long_orbit}
We now describe a Thompson-like flow with very long periodic orbits. 

\begin{minipage}{.40\textwidth}
\begin{center}
\includegraphics[width=.5\textwidth]{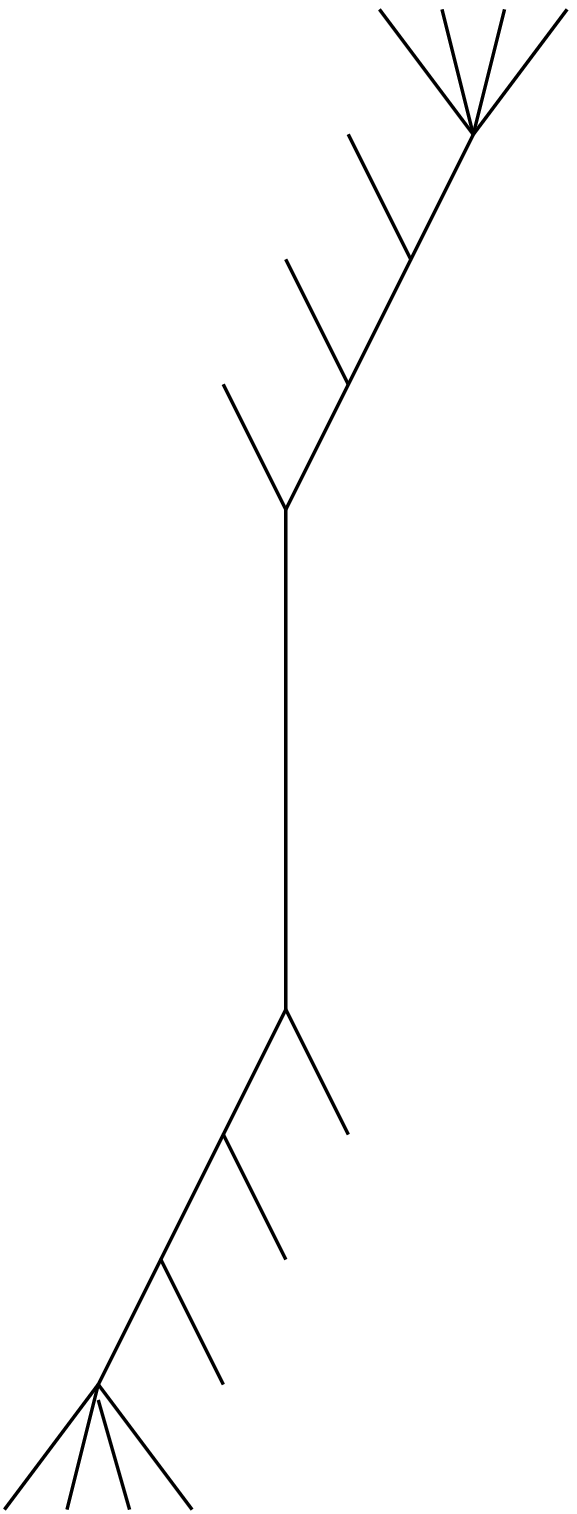}
\relabelbox
\extralabel <-33pt,95pt> {$e$}
\extralabel <-58pt,140pt> {$a_1$}
\extralabel <-50pt,157pt> {$a_2$}
\extralabel <-42pt,174pt> {$a_3$}
\extralabel <0pt,180pt> {$s^+$}
\extralabel <-85pt,10pt> {$s^-$}
\extralabel <-27pt,46pt> {$a_1$}
\extralabel <-35pt,29pt> {$a_2$}
\extralabel <-43pt,12pt> {$a_3$}
\endrelabelbox
\end{center}
\end{minipage}
\begin{minipage}{.56\textwidth}
We define a train track $\tau$ depending on three parameters. A typical example is illustrated in the
figure. There is a middle edge $e$ containing $r_1$ intervals. There are a nested sequence of
$r_2$ switches on either side of $e$, each $3$-valent. The $a_i$ on the top left are glued to the $a_i$
on the bottom right of the figure. Finally, there are two ``bushy" switches $s^\pm$ each of
which is $r_3+1$ valent. The extreme points of these switches are glued with a twist. Note $m = r_1 + 3 r_2 + r_3$.

\vskip 12pt

Each time we split open a sink, the length of the innermost edge doubles. We do this $r_2$ times,
giving an edge of length $r_1 \cdot 2^{r_2}$. When we split open $s^\pm$, the result is a single circle
of length $r_1r_3\cdot 2^{r_2}$. If we set $r_2 = O(m)$ then the length of the periodic orbit is
at least exponential in $m$.
\end{minipage}

\end{example}

The train track in Example~\ref{very_long_orbit} can be embedded in a surface of genus $O(m)$. Its
existence means that one cannot improve the bound in Theorem A' very much without using more
detailed information about the embedding of the train track in a torus.

This suggests some obvious questions:

\begin{question}
Is there a polynomial bound (in $m$) for the denominator of 
rotation number for an element in a (generalized) Thompson's circle group? What about a quadratic bound?
\end{question}

\begin{question}
How does the length of the smallest periodic orbit in a Thompson-like flow
depend on genus?
\end{question}

Finally, the following question seems to be open, and hard to address directly with our methods: 

\begin{question}\label{rational_derivative_question}
Let $r = p/q$ be a non-integral rational number. Let $t$ be an element of $\PL^+(S^1)$ 
mapping rationals to rationals, with rational break points, and with derivatives 
of the form $r^s$ for $s \in \Z$. Is the rotation number of $t$ rational?
\end{question}

\end{document}